\documentclass[11pt, letterpaper]{amsart}
\usepackage{epic, eepic, amsfonts, latexsym,amssymb,graphicx,multicol,mathrsfs,color,xypic,amscd, amsxtra, verbatim, paralist, xspace, url} 
\usepackage{hyperref}
\usepackage[geometry]{ifsym}
\input xy
\xyoption{all}

 \newlength{\baseunit}               
 \newcount{\numlines}                
\numberwithin{equation}{section}

\newenvironment{enumeratea} {\begin{enumerate}[\upshape (a)]} {\end{enumerate}}

\newenvironment{enumerate1} {\begin{enumerate}[\upshape (1)]} {\end{enumerate}}

\newtheorem{lem}{Lemma}[section]
\newtheorem{thm}[lem]{Theorem}

\newtheorem*{thmA}{Theorem A}
\newtheorem*{thmB}{Theorem B}

\newtheorem{prop}[lem]{Proposition}

\newtheorem*{ques*}{Question}
\newtheorem{cor}[lem]{Corollary}

\theoremstyle{definition}
\newtheorem{defn}[lem]{Definition}

\newtheorem{example}[lem]{Example}

\theoremstyle{remark}
\newtheorem{remark}[lem]{Remark}

\newcommand{\epf}{\qed \vspace{+10pt}}

  \newcommand\cC{\mathcal{C}}  \renewcommand\cH{\mathcal{H}}\newcommand\cM{\mathcal{M}}\newcommand\cO{\mathcal{O}}\newcommand\cQ{\mathcal{Q}}\renewcommand\cR{\mathcal{R}}\newcommand\cS{\mathcal{S}}\newcommand\cU{\mathcal{U}}\newcommand\cV{\mathcal{V}}\newcommand\cW{\mathcal{W}}\newcommand\cX{\mathcal{X}}\newcommand\cY{\mathcal{Y}}\newcommand\cZ{\mathcal{Z}}

\renewcommand\AA{\mathbb{A}}\newcommand\CC{\mathbb{C}}\newcommand\GG{\mathbb{G}}\newcommand\PP{\mathbb{P}}
\newcommand\ZZ{\mathbb{Z}}

\newcommand\hookarr{\hookrightarrow}
\newcommand{\xarr}{\xrightarrow}

\renewcommand{\setminus}{\smallsetminus}

\renewcommand{\ss}{\operatorname{ss}}

\newcommand{\rrarrows}{\rightrightarrows}

\newcommand{\Def}{\operatorname{Def}}

\newcommand{\Sym}{\operatorname{Sym}}
\newcommand{\Aut}{\operatorname{Aut}}

\newcommand{\oh}{\cO}

\newcommand{\Spec}{\operatorname{Spec}}

\renewcommand{\tilde}{\widetilde}

\renewcommand{\bar}{\overline}

\newcommand{\PGL}{\operatorname{PGL}}

\newcommand{\dual}{\vee}
\renewcommand{\hat}{\widehat}

\begin{document}

\title{Existence of good moduli spaces for $A_k$-stable curves}

\author[Alper]{Jarod Alper}
\author[Smyth]{David Smyth$^*$}

\address[Alper]{Mathematical Sciences Institute\\
Australia National University\\
Canberra, ACT 0200} 

\address[Smyth]{Department of Mathematics\\
Harvard University\\
1 Oxford Street\\
Cambridge, MA 01238}
\thanks{*The second author was partially supported by NSF grant 
DMS-0901095 during the preparation of this work.}

\begin{abstract}
We prove a general criterion for an algebraic stack to admit a good moduli space. This result may be considered as a weak analog of the Keel-Mori theorem, which guarantees the existence of a coarse moduli space for a separated Deligne-Mumford stack. We apply our result to prove that the moduli stacks of $A_k^-/A_k/A_k^+$-stable curves admit good moduli spaces. In forthcoming work, we will prove that these moduli spaces are projective and use them to construct the second flip in the log minimal model program for $\bar{M}_{g}$.
\end{abstract}

\maketitle

\setcounter{tocdepth}{1}

\section{Introduction}

This is the second in a sequence of three papers, in which we construct the second flip in the log minimal model program for $\bar{M}_{g}$. In this paper, we prove that the moduli stacks  $\bar{\cM}_{g,n}(A_k)$, $\bar{\cM}_{g,n}(A_k^+)$ (introduced in \cite{asw}) admit good moduli spaces  $\bar{M}_{g,n}(A_k)$, $\bar{M}_{g,n}(A_k^+)$ for $k =2,3,4$, and $(g,n) \geq (3,0)$ or $(g,n) \geq (1,1)$.\footnote{Throughout this paper, whenever we write $\bar{\cM}_{g,n}(A_k^-)$, $\bar{\cM}_{g,n}(A_k)$, or $\bar{\cM}_{g,n}(A_k^+)$, we always assume $k =2,3,4$, and $(g,n) \geq (3,0)$ or $(g,n) \geq (1,1)$. } In a forthcoming paper, we will prove that these moduli spaces are projective and interpret them as log canonical models of $\bar{M}_{g,n}$. In the special cases $k=2,3$ and $n=0$, these spaces have been constructed using Geometric Invariant Theory (GIT) by Hassett and Hyeon in \cite{hassett-hyeon_contraction}, \cite{hassett-hyeon_flip}, but there is no known GIT construction of the remaining moduli spaces. Thus, this paper gives the first \emph{intrinsic} construction of a moduli space associated to non-separated stack. The moduli spaces $\bar{M}_{g}(A_4)$ and $\bar{M}_{g}(A_4^+)$ are of particular interest because they constitute the ``target'' and ``flipping space'' respectively of the second flip of the log minimal model program for $\bar{M}_{g}$.  In Section 3, we prove the following result:

\begin{thmA}
$\bar{\cM}_{g,n}(A_k)$ and $\bar{\cM}_{g,n}(A_k^+)$ admit good moduli spaces $\bar{M}_{g,n}(A_k)$ and  $\bar{M}_{g,n}(A_k^+)$ respectively. $\bar{M}_{g,n}(A_k), \bar{M}_{g,n}(A_k^+)$ are proper algebraic spaces.
\end{thmA}

In order to prove Theorem A, we develop an intrinsic technique to prove the existence of good moduli spaces that may be considered as an analog of the Keel-Mori theorem for algebraic stacks. Recall that the Keel-Mori theorem \cite{keel-mori} asserts that \emph{if $\cX$ is a separated Deligne-Mumford stack, then $\cX$ admits a coarse moduli space}. Our main technical result is proven in Section 2:

\begin{thmB}
 Let $\cX$ be an algebraic stack of finite type over $k$.  Suppose that:
\begin{enumerate1}
\item For every closed point $x \in \cX$, there exists an affine, \'etale neighborhood $f: [\Spec A_x / G_x] \to \cX$ of $x$ such that
\begin{enumeratea}
\item $f$ is stabilizer preserving at closed points of $[\Spec A / G_x]$.
\item  $f$ sends closed points to closed points.
\end{enumeratea}
\item For any $k$-point $x \in \cX$, the closed substack $\overline{ \{x\}}$ admits a good moduli space.
\end{enumerate1}
 Then $\cX$ admits a good moduli space. 
 \end{thmB}

\subsection*{Sketch of proof}
In order to motivate the hypotheses of Theorem B as well as our argument, it will be useful to sketch a proof of the Keel-Mori theorem \cite{keel-mori} in the case where $\cX$ is a separated Deligne-Mumford stack of finite type over an algebraically closed field $k$ of characteristic $0$. In this case, every  closed point $x \in \cX$ admits an \'{e}tale neighborhood of the form 
$$[\Spec A_x/G_x] \rightarrow \cX,$$
where $A_x$ is a finite-type $k$-algebra and $G_x$ is the stabilizer of $x$. The union $\coprod_{x \in \cX} [\Spec A_x/G_x] $ defines an etale cover of $\cX$; reducing to a finite subcover, we obtain an atlas $f:\cW \rightarrow \cX$ with the following properties:
\begin{enumerate}
\item $f$ is affine and \'etale,
\item $\cW$ admits a coarse moduli space $W$.
\end{enumerate}
Indeed, (2) follows simply by taking invariants $[\Spec A_x/G_x] \rightarrow \Spec A_x^{G_x}$ and since $f$ is affine, the fiber product $\cR := \cW \times_{\cX} \cW$ admits a coarse moduli space $R$.  
We may thus consider the following diagram:
\begin{equation} \label{diagram-groupoid}
\xymatrix{
\cR  \ar@<.5ex>[r]^{p_1} \ar@<-.5ex>[r]_{p_2} \ar[d]^{\varphi}  & \cW \ar[r]^f \ar[d]^{\phi} \ar[r]  & \cX \\
R  \ar@<.5ex>[r]^{q_1} \ar@<-.5ex>[r]_{q_2}  & W &
}
\end{equation}
The crucial question is: can we choose $f: \cW \to \cX$ to guarantee that the induced projections $q_1, q_2$ are etale? If so, then $R \rightrightarrows W$ defines an \'etale equivalence relation, and the algebraic space quotient $X$ gives a coarse moduli space for $\cX$. The answer, of course, is \emph{yes} because $\cX$ is separated.  Indeed, the condition that $\cX$ is separated implies that the atlas $f$ may be chosen to be stabilizer preserving.\footnote{The set of points $\omega \in \cW$ where $f$ is not stabilizer preserving is simply the image of the complement of the open substack $I_{\cW} \subseteq I_{\cX} \times_{\cX} \cW$ in $\cW$ and therefore is closed since $I_{\cX} \to \cX$ is proper.  By removing this locus from $\cW$, $f: \cW \to \cX$ may be chosen to be stabilizer preserving.} Thus, we may take the projections $\cR \rightrightarrows \cW$ to be stabilizer preserving and \'etale, and this implies that the projections $R \rightrightarrows W$ are \'etale. \footnote{To see this, note that if $r \in R$ is any closed point and $\rho \in \cR$ is its preimage, then $\hat{\oh}_{R,r} \cong D_{\rho}^{G_{\rho}}$, where $D_{\rho}$ denotes the miniversal form deformation space of $\rho$ and $G_{\rho}$ is the stabilizer of $\rho$; similarly $\hat{\oh}_{W, q_i(r)} \cong D_{p_i(\rho)}^{G_{p_i(\rho)}}$.  Now $p_i$ \'etale implies $D_{\rho} \cong D_{p_i(\rho)}$ and $p_i$ stabilizer preserving implies $G_{\rho} \cong G_{p_i(\rho)}$, so $\hat{\oh}_{R,r} \cong \hat{\oh}_{W, q_i(r)}$, i.e. $q_i$ is \'etale.}

Now let us ask whether this proof can be adapted to construct a good moduli space when $\cX$ is an Artin stack of finite type over an algebraically closed field $k$. The first complication that arises is the following:

\begin{ques*}
 Does any closed point $x \in \cX$ admit an etale neighborhood of the form $[\Spec A_x/G_x] \rightarrow \cX$
 \end{ques*}

This is an open problem  \cite[Conjecture 1]{alper_quotient}, but if we make the mild assumption that $\cX$ is a global quotient stack of a normal scheme by a connected algebraic group with reductive stabilizers at closed points, then the desired \'etale neighborhoods exist (Proposition \ref{prop-quotient}). Precisely as in the case of the Deligne-Mumford stack, we thus obtain a cover $\cW \rightarrow \cX$ satisfying properties $(1)$-$(2)$ above, except that we must replace ``coarse moduli space'' by ``good moduli space.'' Exactly as in the case of the Keel-Mori theorem, if the projections $q_1,q_2$ define an etale equivalence relation, the resulting algebraic space quotient $X$ will be a good moduli space for $\cX$. The natural question is therefore: What hypotheses must be imposed on $\cX$ or the presentation $f:\cW \rightarrow \cX$ to ensure an \'{e}tale equivalence relation? More precisely, can we identify a sufficient set of hypotheses which can be directly verified for geometrically-defined stacks, such as $\bar{\cM}_{g,n}(A_k)$ and $\bar{\cM}_{g,n}(A_k^+)$?

Our result gives at least one plausible answer to these questions.   With notation as in Diagram  \ref{diagram-groupoid}, where $\cW \to W$ and $\cR \to R$ are now good moduli spaces, let us examine what conditions are necessary to conclude that $q_1, q_2$ are \'etale.  
One difficulty is that $f: \cW \rightarrow \cX$ need not be stabilizer-preserving, since $\cX$ is a non-separated algebraic stack. A second difficulty is that if $\omega \in \cW$ is a $k$-point with image $w \in W$, then the formal neighborhood $\hat{\oh}_{W,w}$ can be identified with the invariants $D_{\omega}^{G_{\omega}}$ of the miniversal deformation space if and only if $\omega \in \cW$ is a \emph{closed point}.  Thus, if $r \in R$ is any $k$-point, then for $q_1, q_2$ to be \'etale at $r$, or equivalently for the induced maps $\hat{\oh}_{W,w} \to \hat{\oh}_{R,r}$ to be isomorphisms, we must manually impose the following conditions:  $p_1, p_2$ should be stabilizer preserving at $\rho$ and $p_1(\rho)$, $p_2(\rho)$ should be closed points, where  $\rho \in \cR$ is the unique closed point in the preimage of $r \in R$.  We have now identified two key conditions that will imply that $R \rightrightarrows W$ is an \'etale equivalence relation:
\begin{enumerate}
\item[($\star$)] The morphism $f: \cW \to \cX$ is stabilizer preserving at closed points.
\item[($\star \star$)] The projections $p_1, p_2: \cW \times_{\cX} \cW$ send closed points to closed points.
\end{enumerate}

Condition ($\star$) is precisely hypothesis (2a) of Theorem B.  In practice, it is difficult to directly verify condition ($\star \star$), but it turns out that it is implied by conditions (2b) and (3), which are often easier to verify.

Our main application of Theorem B is to construct the good moduli spaces for the moduli stacks of curves $\bar{\cM}_{g,n}(A_k), \bar{\cM}_{g,n}(A_k^+)$ (Theorem A). Indeed, using our explicit classification of closed points of $\bar{\cM}_{g,n}(A_k), \bar{\cM}_{g,n}(A_k^+)$ and the \'etale local VGIT description of the inclusions $\bar{\cM}_{g,n}(A_k^-) \hookrightarrow \bar{\cM}_{g,n}(A_k) \hookleftarrow \bar{\cM}_{g,n}(A_k^+)$ \cite[Section 7]{asw}, we can directly verify hypotheses (1a), (1b), and (2) of Theorem B. The interested reader who wishes to understand the \emph{geometric} significance of these hypotheses is urged to read Appendix A, where we have included three examples of moduli stacks of curves which fail to satisfy (1a), (1b) and (2) respectively. Studying these examples should convince the reader that if one thinks of a moduli stack of curves $\cM$ as given by ``crushing'' a certain subcurves $Z(C)$ of each stable curve $C$, then hypotheses (1a), (1b) express (very roughly) the following conditions,\footnote{Significantly, these are precisely the defining conditions of the extremal assignments occurring in the combinatorial classification of stable modular classifications \cite{smyth_zstable}.} which are easily seen to hold for the subcurves that we replace in our definitions of $\bar{\cM}_{g,n}(A_k), \bar{\cM}_{g,n}(A_k^+)$.
\begin{itemize}
\item[(1a)] The subcurve $Z(C)$ should be $\Aut(C)$ invariant.
\item[(1b)] If $C \leadsto C'$ is a specialization of stable curves, then a component of $C$ is in $Z(C)$ if and only if its closure is contained in $Z(C')$. 
\end{itemize}
Condition (2) is rather different; roughly speaking, it asserts that for any $k$-point $x \in \cX$, the closed substack $\overline{ \{x\}}$ should look like the quotient of an affine scheme by a reductive group.


\subsection*{Notation}  We assume throughout that $k$ is an algebraically closed field of characteristic $0$ and that all algebraic stacks are quasi-separated. A curve is a reduced, connected, 1-dimensional scheme of finite type over $k$. $\Delta$ will always denote the spectrum of a valuation ring $R$ with field of fractions $K$, and we let $0 \in \Delta$ denote the closed point and $\eta \in \Delta$ the generic point.  We will also use the following definitions:

\begin{defn} Let $f: \cX \to \cY$ be a morphism of algebraic stacks of finite type over $k$.  We say that
\begin{itemize}
\item $f$ \emph{sends closed points to closed points} if for every closed point $x \in \cX$, $f(x) \in \cY$ is closed.
\item $f$ is \emph{stabilizer preserving at $x \in \cX(k)$} if $\Aut_{\cX(k)}(x) \to \Aut_{\cY(k)}(f(x))$ is an isomorphism.
\item $f$ is \emph{stabilizer preserving} if $I_{\cX} \to I_{\cY} \times_{\cY} \cX$ is an isomorphism.
\end{itemize}
\end{defn}

\begin{defn}
If $\phi: \cX \to X$ is a good moduli space, we say that an open substack $\cU \subseteq \cX$ is \emph{saturated} if $\phi^{-1}(\phi(\cU)) = \cU$.  This is equivalent to requiring that the open immersion $\cU \hookarr \cX$ sends closed points to closed points.
\end{defn}

\begin{defn} \label{definition-etale-presentations}  Let $\cX$ be an algebraic stack of finite type over $k$ and $x \in \cX(k)$.   We call
$f: (\cW, w) \to \cX$ an \emph{\'etale GIT presentation around $x$} if 
 \begin{itemize}
 \item $\cW = [\Spec A /G_x]$ where $A$ is a finite type $k$-algebra and $G_x$ is linearly reductive.
 \item $f: \cW \to \cX$ is an \'etale, affine morphism with $f(w) = x$.
 \item $f$ is stabilizer preserving at $w \in \cW$.
 \end{itemize}
We say that $\cX$ \emph{admits \'etale GIT presentations} if there exist \'etale GIT presentations around all closed points $x \in \cX(k)$. 
\end{defn}

\section{Proof of Theorem B}

\subsection*{Existence of \'etale GIT presentations}
We begin by including the following application of Sumihiro's Theorem \cite{sumihiro1} and Luna's \'Etale Slice Theorem \cite{luna}.  This argument appeared in \cite[Theorem 3]{alper_quotient} but given it's importance to the main results of this paper and it's simplicity, we include the proof.

\begin{prop} \label{prop-quotient}
 Let $\cX$ be an algebraic stack finite type over $k$, and suppose that $\cX$ is a quotient stack $[X/G]$ where $G$ is a connected algebraic group acting on a normal separated scheme $X$.  If $x \in \cX(k)$ has linearly reductive stabilizer, there exists a locally closed $G_x$-invariant affine $W \hookarr X$ with $w \in W$ such that
$$[W/G_x] \to [X/G]$$
is affine and \'etale.  In particular, if $\cX$ has linearly reductive stabilizers at all closed points $x \in \cX(k)$, then $\cX$ admits \'etale GIT presentations.
\end{prop}

\begin{proof}
By applying \cite[Theorem 1 and Lemma 8]{sumihiro1}, there exists an open $G$-invariant $U_1$ containing $x$ and a $G$-equivariant immersion $U_1 \hookarr Y= \PP(V)$ where $V$ is a $G$-representation.  Since the action of $G_x$ on $\Spec \Sym^* V^{\dual}$ fixes the line spanned by $x$, there exists a $G_x$-invariant homogeneous polynomial $f$ with $f(x) \neq 0$.  Then $Y_f$ is a $G_x$-invariant affine.  $G_x$ acts on $T_x Y_{f}$ and there exists a $G_x$-invariant morphism $g: Y_f \to T_x \PP(V)$ which is \'etale.   Since $G_x$ is linearly reductive, we may write $T_x \PP(V) = T_x o(x) \oplus N$ for a $G_x$-representation $N$.  Then $W_1 = g^{-1}(N) \cap U_1$ is a $G_x$-invariant closed quasi-affine subscheme and we may choose a $G_x$-invariant affine $W \subseteq W_1$ containing $x$.  It is easy to see that $[W/G_x] \to [X/G]$ is \'etale at $w$.  By shrinking $W$ further, one obtains that $[W/G_x] \to [X/G]$ is \'etale which is clearly affine.
 \end{proof}

\subsection*{\'Etale descent}
We will need the following characterization of when an \'etale morphism of algebraic stacks induces an \'etale morphism of good moduli spaces.

\begin{thm} (\cite[Theorem 5.1]{alper_good}) \label{etale-preserving}
Consider a commutative diagram $$\xymatrix{ 
\cW \ar[r]^f \ar[d]^{\varphi}		& \cX \ar[d]^{\phi} \\
W \ar[r]^g					& X
}$$
where $f$ is a representable morphism between algebraic stacks of finite type over $k$ and $\varphi, \phi$ are good moduli spaces.  Let $w \in |\cW|$.  Suppose, $f$ is \'etale and stabilizer preserving at $w$, and that both $w$ and $f(w)$ are closed.
Then $g$ is \'etale at $\varphi(w)$. \epf
\end{thm}

\begin{cor} \label{corollary-cartesian}
Consider a commutative diagram $$\xymatrix{ 
\cW \ar[r]^f \ar[d]^{\varphi}		& \cX \ar[d]^{\phi} \\
W \ar[r]^g					& X
}$$
where $f$ is a morphism between algebraic stacks of finite type over $k$ and $\varphi, \phi$ are good moduli spaces.  Suppose $f$ is representable, separated, \'etale, sends closed points to closed points, and is stabilizer preserving at closed points of $|\cW|$.  Then $g$ is \'etale and the above diagram is cartesian.
\end{cor}

\begin{proof}
Theorem \ref{etale-preserving} implies that $g$ is \'etale.  The hypotheses imply that the induced morphism $\Psi: W \times_X \cX \to \cX$ is representable, separated, quasi-finite and sends closed points to closed points.  Therefore \cite[Proposition 6.4]{alper_good}) implies that $\Psi$ is a finite.  Moreover, since $f$ and $g$ are \'etale, so is $\Psi$.  But since $\cW$ and $W \times_X \cX$ both have $W$ as a good moduli space, it follows that a closed point in $W \times_X \cX$ has a unique preimage under $\Psi$.  Therefore, $\Psi$ is an isomorphism and the diagram is cartesian.
\end{proof}

\begin{cor} \label{corollary-openness}
Consider a commutative diagram 
$$\xymatrix{ 
\cW \ar[r]^f \ar[d]^{\varphi}		& \cX \ar[d]^{\phi} \\
W \ar[r]^g					& X
}$$
where $f$ is a representable, separated and \'etale morphism between algebraic stacks of finite type over $k$ and $\varphi, \phi$ are good moduli spaces.  Then there exists a saturated open substack $\cU \subseteq \cW$ such that
\begin{enumeratea}
\item The induced diagram
$$\xymatrix{ 
\cU \ar[r]^{f|_{\cU}} \ar[d]^{\varphi|_{\cU}}		& \cX \ar[d]^{\phi} \\
\varphi(\cU) \ar[r]^{g|_{\varphi(\cU)}}					& X
}$$
is cartesian and $g|_{\varphi(\cU)}$ is \'etale.  In particular, $f|_{\cU}$ is stabilizer preserving and sends closed points to closed points.
\item If $w \in \cW$ is a closed point such that $f$ is stabilizer preserving at $w$ and $f(w) \in \cX$ is closed, then $w \in \cU$.
\end{enumeratea}
\end{cor}

\begin{proof}  By Zariski's Main Theorem, we may factor $f$ as 
$$f: \cW \xarr{i} \tilde{\cW} \xarr{\tilde{f}} \cX$$
where $i$ is an open immersion and $\tilde f$ is finite.  Since $\cX$ admits a good moduli space and $\tilde f: \tilde{\cW} \to \cX$ is finite, there exists a good moduli space $\tilde{\varphi}: \tilde{\cW} \to \tilde{W}$.  Let $\cZ = \tilde{\varphi}^{-1}(\tilde{\varphi}(\tilde{\cW} \setminus \cW))$ and $\cU' = \tilde{\cW} \setminus \cZ \subseteq \cW$.   Observe that a point $w \in \cW$ is in $\cZ$ if and only if $\overline{ \{w \} } \cap (\tilde{\cW} \setminus \cW) \neq \emptyset$.  If $u \in \cU'$ is a closed point, then $u \in \tilde{\cW}$ is also closed and since  $f: \tilde{\cW} \to \cX$ is finite, $f(u) \in \cX$ is closed.  
Therefore, $f|_{\cU'}$ sends closed points to closed points.  
Observe first that if $f$ is stabilizer preserving at a closed point $u \in \cU'$, then by Theorem \ref{etale-preserving}, $g$ is \'etale at $\varphi(u)$.   If we set $V \subseteq \phi(\cU')$ is the open locus where $g|_{V}$ is \'etale, then $\varphi^{-1}(V)$ contains $u_0$ and therefore $u$.  Thus, we can shrink $\cU'$ so that $g|_{\varphi(\cU')}$ is \'etale.  This gives a diagram
$$\xymatrix{ 
\cU' \ar[r]^{f|_{\cU'}} \ar[d]^{\varphi|_{\cU'}}		& \cX \ar[d]^{\phi} \\
\varphi(\cU') \ar[r]^{g|_{\varphi(\cU')}}					& X
}$$
The induced map $\Psi: \cU' \to \varphi(\cU') \times_X \cX$ is affine and \'etale, and also finite by \cite[Proposition 6.4]{alper_good}).  If $f$ is stabilizer preserving at $u \in \cU'$, then $\Psi$ has degree $1$ at $u$ so that there is an open locus $\cU \subseteq \cU'$ where $\Psi$ is an isomorphism.  This establishes (a) as well as (b).
\end{proof}
The following simple lemma will allows us to shrink Zariski-locally by saturated open substacks.

\begin{lem} \label{lemma-shrinking}
Let $\phi: \cX \to X$ be a good moduli space.  Let $x \in \cX$ be a closed point and $\cU \subseteq \cX$ be an open substack containing $x$.  Then exists a saturated open substack $\cU_1 \subseteq \cU$ containing $x$.
\end{lem}

\begin{proof}  The substacks $\{x\}$ and $\cX \setminus \cU$ are closed and disjoint.  By \cite[Theorem 4.16]{alper_good}, $\phi(\{x\})$ and $Z:=\phi(\cX \setminus \cU)$ are closed and disjoint.   Therefore, we take $\cU_1 = \phi^{-1}(X \setminus Z)$.
\end{proof}

\subsection*{General existence results}

We begin with a proposition giving conditions on when the good moduli spaces of \'etale GIT presentations can be glued.  A weaker version of this result appear in \cite{alper_local}. 

\begin{prop} \label{prop-existence}
Let $\cX$ be an algebraic stack locally of finite type over $k$.  Suppose there exists an affine, \'etale and surjective morphism $f: \cW \to \cX$ from an algebraic stack $\cW$ such that:
\begin{enumerate1}
\item $\cW$ admits a good moduli space $\varphi: \cW \to W$.
\item $f$ is stabilizer preserving at closed points in $\cW$.
\item The projections $\cW \times_{\cX} \cW \rightrightarrows \cW$ sends closed points to closed points.
\end{enumerate1}
Then there exists a good moduli space $\phi: \cX \to X$ inducing a cartesian diagram
$$\xymatrix{
\cW \ar[r]^f \ar[d]^{\varphi}   & \cX \ar@{-->}[d]^{\phi} \\
W \ar@{-->}[r]  ^g              & X
}$$
such that $g: W \to X$ is \'etale.
\end{prop}

\begin{proof} 
Set $\cX_1 = \cW$, $X_1 = W$  and $\cX_2 = \cX_1 \times_{\cX} \cX_1$ with projections $p_1, p_2:  \cX_2 \to \cX_1$. 
Since $f$ is affine, there exists a good moduli space $\phi_2: \cX_2 \to X_2$.  The two projections $p_1, p_2: \cX_2 \to \cX_1$ induce two morphisms $q_1, q_2: X_2 \to X_1$ such that $q_i \circ \phi_2 = \phi_1 \circ p_i$ for $i=1,2$.  By Theorem \ref{etale-preserving} and hypotheses (2) and (3), $q_1$ and $q_2$ are \'etale.  The induced morphisms $\cX_2 \to X_2 \times_{q_i, X_1, \phi_1} \cX_1$ are isomorphisms by Corollary \ref{corollary-cartesian}.  
Similarly, by setting $\cX_3 = \cX_1 \times_{\cX} \cX_1 \times_{\cX} \cX_1$, there is a good moduli space $\phi_3: \cX_3 \to X_3$ where $X_3 = X_2 \times_{q_1, X_1, q_2} X_2$ and an induced diagram
$$\xymatrix{
\cX_3  \ar@<1ex>[r] \ar@<-1ex>[r] \ar[r] \ar[d]^{\phi_3}    & \cX_2 \ar@<.5ex>[r] \ar@<-.5ex>[r] \ar[d]^{\phi_2}    & \cX_1 \ar[r]^f \ar[d]^{\phi_1}    & \cX \\
X_3  \ar@<1ex>[r] \ar@<-1ex>[r] \ar[r]  & X_2 \ar@<.5ex>[r] \ar@<-.5ex>[r]  & X_1
}$$
where the appropriate squares are cartesian.  Moreover, by the universality of good moduli spaces, there is an induced identity map $X_1 \to X_2$, an inverse $X_2 \to X_2$ and a composition $X_2 \times_{q_1, X_1, q_2} X_2 \to X_2$  giving $X_2 \rrarrows X_1$ an \'etale groupoid structure.

To check that $\Delta: X_2 \to X_1 \times X_1$ is a monomorphism, it suffices to check that there is a unique pre-image of $(x_1, x_1) \in X_1 \times X_1$ where $x_1 \in X_1(k)$.  Let $
\xi_1 \in |\cX_1|$ be the unique closed point in $\phi_1^{-1}(x_1)$.  Since $\cX_1 \to \cX$ is stabilizer preserving at $\xi_1$, we can set $G := \Aut_{\cX_1(k)}(\xi_1) \cong \Aut_{\cX(k)}(f(\xi_1))$.  There are diagrams
$$\xymatrix{
BG \ar[r] \ar[d]                & BG \times BG \ar[d] \\
\cX_2 \ar[r] \ar[d]         & \cX_1 \times \cX_1 \ar[d] \ar@{}[ul] |{\square} \\
\cX \ar[r]                  & \cX \times \cX \ar@{}[ul] |{\square}
}
\qquad
\qquad
\xymatrix{
\cX_2 \ar[r]^{(p_1,p_2)} \ar[d]^{\phi_2}	& \cX_1 \times \cX_1 \ar[d]^{\phi_1 \times \phi_1}\\
X_2 \ar[r]^{\Delta}			& X_1 \times X_1
}
$$
where the squares in the left diagram are 2-cartesian.  Suppose $x_2 \in X_2(k)$ is in the preimage of $(x_1, x_1)$ under $\Delta: X_2 \to X_1 \times X_1$.  Let $\xi_2 \in |\cX_2|$ be the unique closed point in $\phi_{2}^{-1}(x_2)$.  Then $(p_1(\xi_2), p_2(\xi_2)) \in |\cX_1 \times \cX_1|$ is closed and is therefore the unique closed point $(\xi_1, \xi_1)$ in the $(\phi_1 \times \phi_1)^{-1}(x_1,x_1)$.
But since the above diagram is cartesian, $\xi_2$ is the unique point in $|\cX_2|$ which maps to $(\xi_1, \xi_1)$ under  $\cX_2 \to \cX_1 \times \cX_1$.  Therefore, $x_2$ is the unique preimage of $(x_1,x_1)$.

Since $X_2 \times_{q_1, X_1, q_2} X_2 \to X_2$ is an \'etale equivalence relation, there exists an algebraic space quotient $X$ and induced maps $\phi: \cX \to X$ and $X_1 \to X$.  Consider
$$\xymatrix{
\cX_2 \ar[r] \ar[d]  & \cX_1 \ar[r] \ar[d]     & X_1 \ar[d] \\
\cX_1 \ar[r] &\cX \ar[r]              & X
}$$
Since $\cX_2 \cong \cX_1 \times_{X_1} X_2$ and $X_2 \cong X_1 \times_X X_1$, the left and outer square above are 2-cartesian.  Since $\cX_1 \to \cX$ is \'etale and surjective, it follows that the right square is cartesian.  By descent, $\phi: \cX \to X$ is a good moduli space.   \end{proof}

\begin{thm} \label{theorem-global-existence}
 Let $\cX$ be an algebraic stack of finite type over $k$.  Suppose that:
\begin{enumerate1}
\item For every closed point $x \in \cX$, there exists an \'etale GIT presentation $f: \cW \to \cX$ around $x$ such that
\begin{enumeratea}
\item $f$ is stabilizer preserving at closed points in $\cW$.
\item  $f$ sends closed points to closed points
\end{enumeratea}
\item For any $k$-point $x \in \cX$, the closed substack $\overline{ \{x\}}$ admits a good moduli space.
\end{enumerate1}
 Then there exists a good moduli space $\cX \to X$.
\end{thm}

\begin{proof}
Choose \'etale GIT presentations $[\Spec A_x / G_x] \to \cX$ around each closed $x \in \cX$.  Then by refining  the \'etale presentation $\coprod_{x \in \cX} [\Spec A_x/G_x] \to \cX$ to a finite subcover, we obtain an \'etale, affine and surjective morphism $f: \cW \to \cX$ where $\cW$ admits a good moduli space and $f$ is both stabilizer preserving at closed points in $\cW$ and sends closed points to closed points.  We now check that property (2) assures that the projections $p_1, p_2: \cR:=\cW \times_{\cX} \cW \to \cW$ send closed points to closed points. Once this is accomplished, the statement follows from Proposition \ref{prop-existence}.

Let $\rho \in \cR$ be a closed point and set $x = f(p_1(\rho)) = f(p_2(\rho)) \in \cX$.  Let $\cZ = \overline{ \{x\}} \subseteq \cX$, $\cW' = \cW \times_{\cX} \cZ$ and $\cR' = \cR \times_{\cX} \cZ$ with induced maps $f': \cW' \to \cZ$ and $p_1', p_2': \cR' \to \cW'$.  Consider
$$\xymatrix{
\cR'  \ar@<.5ex>[r]^{p_1'} \ar@<-.5ex>[r] _{p_2'} \ar@{^(->}[d]	& \cW' \ar[r]^{f'} \ar@{^(->}[d]	& \cZ \ar@{^(->}[d] \\
\cR  \ar@<.5ex>[r]^{p_1} \ar@<-.5ex>[r]_{p_2} 			& \cW \ar[r]^{f}			& \cX 
}$$
The morphism $f': \cW' \to \cX'$ is \'etale, surjective, affine, sends closed points to closed points and is stabilizer preserving at closed points.  Moreover, by hypothesis (2), there is a diagram
$$\xymatrix{
\cW' \ar[d]^{\phi'} \ar[r]	& \cZ \ar[d]^{\varphi'}  \\
W' \ar[r]^{g'}			& \Spec k
}$$
where $\phi': \cW' \to W'$ and $\varphi': \cZ' \to \Spec k$ are good moduli spaces.  By Corollary \ref{corollary-cartesian}, the above diagram is cartesian.  It follows that the projects $p_1', p_2': \cR' \rrarrows \cW'$ send closed points to closed points as they are the base change of the projections $W' \times W' \to W'$ by $ \cW' \to W'$.  Therefore $p_1'(\rho), p_2'(\rho) \in \cW'$ are closed or equivalently $p_1(\rho), p_2(\rho) \in \cW$ are closed. \end{proof}

\section{Proof of Theorem A}

\subsection*{General facts}

The following lemmas will  under suitable hypotheses under which the existence of a finite surjective morphism $\cX \to \cY$ from an algebraic stack $\cX$ admitting a good moduli space implies that $\cY$ admits a good moduli space.

\begin{lem} \label{lemma-proper-existence}
Let $f: \cX \to \cY$ be a morphism of algebraic stacks of finite type over $k$.  Suppose that:
\begin{enumerate1}
\item $f: \cX \to \cY$ is finite, surjective, and stabilizer preserving at closed points.
\item $\cX$ admits a good moduli space.
\item  $\cY$ admits \'etale GIT presentations.
\end{enumerate1}
Then $\cY$ admits a good moduli space.
 \end{lem}
 \begin{proof}  Let $y \in \cY$ be a closed point and choose an \'etale GIT presentation $g: ([\Spec A_y / G_y],y') \to \cY$ around $y$.  Consider the 2-cartesian diagram
 $$\xymatrix{
 \cX' \ar[r]^{f'} \ar[d]^{g'}       & [\Spec A_y / G_y] \ar[d]^g \\
 \cX \ar[r]^f                   & \cY \ar@{}[ul] |{\square}
 }$$
Since $g$ sends closed points to closed points and $f$ and $f'$ are finite, $g'$ sends closed points to closed points.
 Since $g$ is affine, $\cX'$ admits a good moduli space.
 We note also that $g'$ is stabilizer preserving at every $x' \in f'^{-1}(y')$.  By Corollary \ref{corollary-openness}, 
 there exists an open substack $\cU' \subseteq \cX'$ containing $f'^{-1}(y')$ such that $g'|_{\cU'}$ is stabilizer preserving and sends closed points to closed points.
  Then $\cV' = \cX' \setminus f'(\cX' \setminus \cU')$ is an open substack containing $y'$ and by Lemma \ref{lemma-shrinking} we may shrink $\cV'$ so that it is a saturated open substack containing $y'$.  The pre image $f'^{-1}(\cV')$ is a saturated open substack of $\cX'$ and since $f'^{-1}(\cV') \to \cX$ is the base change of a morphism of algebraic spaces, the projections $f'^{-1}(\cV') \times_{\cX} f'^{-1}(\cV') \rrarrows f'^{-1}(\cV')$ send closed points to closed points.  Therefore, we may assume that $g': \cX' \to \cX$ is stabilizer preserving and that the projections $\cX' \times_{\cX} \cX'$ send closed points to closed points.
By Proposition \ref{prop-existence}, we conclude that $g([\Spec A_y / G_y])$ admits a good moduli space.  But since $f^{-1}(g([\Spec A_y / G_y]))$ is a saturated open substack of $\cX$, it follows that $g([\Spec A_y / G_y]) \hookarr \cY$ sends closed points to closed points.  It follows from \cite[Proposition 7.9]{alper_good} that good moduli spaces for $g([\Spec A_y / G_y])$ may be glued to construct a good moduli space for $\cY$.
 \end{proof}

 \begin{lem} \label{lemma-finite-affine-existence}
 Consider a commutative diagram
 $$\xymatrix{
 \cX \ar[r] \ar[rd] & \cY \ar[d] \\
                & Z
}$$
of algebraic stacks of finite type over $k$ where $Z$ is an algebraic space.  Suppose that:
\begin{enumerate1}
\item $\cX \to \cY$ is finite and surjective.
\item $\cX \to Z$ is cohomologically affine.
\item $\cY$ is a quotient stack of the form  $[Y/G]$ where $Y$ is an algebraic space with an action of a linearly reductive group $G$.) 
\end{enumerate1}
Then $\cY \to Z$ is cohomologically affine.
 \end{lem}

 \begin{proof}  Since $\cX \to \cY$ is affine, $\cX$ is the quotient stack $\cX = [X/G]$ where  $X = \cX \times_{\cY} Y$.  Since $X \to \cX$ is affine and $\cX \to Z$ is cohomologically affine, $X \to Z$ is affine by Serre's criterion.  The morphism $X \to Y$ is a finite and surjective so by Chevalley's theorem, we can conclude that $Y \to Z$ is affine.  Therefore $\cY \to Z$ is cohomologically affine.
 \end{proof}

\begin{prop} \label{prop-finite-existence}
Let $f: \cX \to \cY$ be a morphism of algebraic stacks of finite type over  $k$.  Suppose that:
\begin{enumerate1}
\item $\cX \to \cY$ is finite and surjective.
\item There exists a good moduli space $\cX \to X$ with $X$ separated.
\item  $\cY$ is a quotient stack and admits \'etale GIT presentations.
\end{enumerate1}
Then $\cY$ admits a good moduli space.
\end{prop}

\begin{proof} Let $\cZ$ be the scheme-theoretic image of $\cX \to X \times \cY$.  Since $\cX \to \cY$ is finite, so is $\cX \to \cZ$.  By Lemma \ref{lemma-finite-affine-existence}, the projection $\cZ \to X$ is cohomologically affine so that $\cZ$ admits a good moduli space.  The composition $\cZ \hookarr X \times \cY \to \cY$ is finite, surjective and stabilizer preserving.  
The statement now follows from Lemma \ref{lemma-proper-existence}.
\end{proof}

\begin{remark}
They hypothesis (2) that $X$ is separated is necessary.  For example, let $X$ be the affine line with the point $0$ doubled and let $\ZZ_2$ act on $X$ by swapping the points at $0$ and fixing all other points.  Then $X \to [X/\ZZ_2]$ satisfies the hypotheses but $[X/\ZZ_2]$ does not admit a good moduli space.
\end{remark}

\subsection*{Inductive existence}
The following lemma will be used to show that if $\cX^- = \bar{\cM}_{g,n}(A_k^-)$, $\cV^- = \bar{\cS}_{g,n}(A_k) = \bar \cM_{g,n}(A_k) \setminus \bar{\cM}_{g,n} (A_k^-)$  and $\cV^+ = \bar{\cH}_{g,n}(A_k) = \bar \cM_{g,n}(A_k) \setminus \bar{\cM}_{g,n} (A_k^+)$ admit good moduli spaces, then so does $\cX=\bar{\cM}_{g,n}(A_k)$.

\begin{prop} \label{prop-inductive-existence1}
Let $\cX$ be an algebraic stack of finite type over $k$.  Let $\cX^- \subseteq \cX$ be an open substack and $\cV^-,  \cV^+ \subseteq \cX^-$ be closed substacks where $\cX^-$ is the complement of $\cV^-$.  Suppose:
\begin{enumerate1}
\item There exist good moduli spaces $\cX^- \to X^-$, $\cV^- \to V^-$ and $\cV^{+} \to V^+$.
\item \label{part-closed}
For every closed point $z \in |\cX \setminus \cX^-|$, $z \in \cV^+$.
\item \label{part-closed2}
For every closed point $z \in |\cX^-|$ which is not closed in $|\cX|$, $z \in \cV^+$.
\item \label{part-stable-reduction} If $x \in \cX \setminus \cX^-$, then for every diagram
$$
\xymatrix{
 \Spec K \ar[r] \ar[d]  & \cX^{-} \ar@{^(->}[d] \\
 \Spec R \ar[r]^h \ar[ur]^{\tilde h}  & \cX
}$$
such that $h(0) = x$ and $\tilde h(0) \in |\cX^{-}|$ closed, we have $\tilde h(0) \in \cV^+$.
\item For any $k$-point $x \in \cX$, the closed substack $\overline{ \{x\}}$ admits a good moduli space.
\item For any closed point $x \in \cX$ there exists an \'etale GIT presentation $f:(\cW,w) \rightarrow \cX$ around $x \in \cX$ inducing a diagram
$$\xymatrix{
\cW^{-} \ar[d]_{f|_{\cW^-}} \ar@{^(->}[r]      & \cW \ar[d]^f \\
\cX^{-} \ar@{^(->}[r]^{\text{op}}                   & \cX       & \ar@{_(->}[l]_{\text{cl}} \cV^{-}, \cV^{+}
}$$
where $\cW^{-} = f^{-1}(\cX^-)$ such that there exist good moduli spaces $\cW^- \to W^-$ and $\cW \to W$, and $W^- \to W$ is proper.

\end{enumerate1}
Then there exists a good moduli space $\cX \to X$ such that $X^- \to X$ is proper.
\end{prop}

\begin{proof}
Let $x \in \cX$ be a closed point and $f:(\cW,w) \rightarrow \cX$ be an \'etale GIT presentation satisfying hypothesis (6).  We claim that, without loss of generality, we may assume that
\begin{enumerate}
\item $f|_{f^{-1}(\cV^-)}$, $f|_{f^{-1}(\cV^+)}$ is stabilizer preserving and sends closed points to closed points.
\item $f|_{\cW^-}$ is stabilizer preserving and sends closed points to closed points.
\end{enumerate}

For (1), note that since $\cV^-$ admits a good moduli space, Corollary \ref{corollary-openness} implies the existence of an open substack $\cQ \subseteq f^{-1}(\cV^{-})$ containing $w$ such that $f|_{\cQ}$ is stabilizer preserving and sends closed points to closed points.  Let $\cU \subseteq \cW$ be an open substack such that $\cU \cap f^{-1}(\cV^-) = \cQ$.  By Lemma \ref{lemma-shrinking}, we may shrink $\cU$ to a saturated open substack $\cU'$ of $\cW$ containing $w$. Then $\cU' \cap f^{-1}(\cV^-) \subseteq f^{-1}(\cV^-)$ and $\cU' \cap \cW^- \subseteq \cW^-$ are saturated open substacks.  By replacing $\cW$ with $\cU'$, we may assume that $f|_{f^{-1}(\cV^-)}$ is stabilizer preserving and sends closed points to closed points.  Similarly, we may assume that $f|_{f^{-1}(\cV^+)}$ is stabilizer preserving and sends closed points to closed points.

For (2), by Corollary \ref{corollary-openness}, there exists an open substack $\cU \subseteq \cW^-$ such that $f|_{\cU}: \cU \to \cX^-$ is stabilizer preserving and sends closed points to closed points; moreover, $\cU$ contains all closed points $w \in \cW^-$ such that $f$ is stabilizer preserving at $w$ and $f(w) \in \cX^-$ is closed.  Let $\cZ = \cW^{-} \setminus \cU$ and $\bar{\cZ}$ be the closure of $\cZ$ in $\cW$.  We claim that $w \notin \bar{\cZ}$. Once this is established, Lemma \ref{lemma-shrinking} implies that we may replace shrink $\cW$ by a saturated open substack containing $w$ so that $f|_{\cW^-}$ is stabilizer preserving and sends closed points to closed points.  If $x \in \cX^-$, then clearly $w \notin \bar{\cZ}$ so we may suppose that $x \in \cX \setminus \cX^-$.  Suppose, by way of contradiction, that $w \in \bar{\cZ}$. Then there exists a specialization diagram
$$\xymatrix{
\Spec K \ar[r] \ar[d]       & \cZ \ar[d] \\
\Spec R \ar[r]^h            & \cW
}$$
such that $h(0) = w$. Since the composition $\cW^{-} \to W^- \to W$ is universally closed, after an extension of the fraction field $K$, there exists a diagram
$$
\xymatrix{
 \Spec K \ar[r] \ar[d]          & \cW^{-} \ar@{^(->}[d] \ar[r]      & W^{-} \ar[d]  \\
 \Spec R \ar[r]^h \ar@{-->}[ur]^{\tilde h}  & \cW \ar[r]    & W
}$$
where the lift $\tilde h: \Spec R \to \cW^{-}$ extends $\Spec K \to \cW^{-}$ and $\tilde h(0) \in \cW^-$ is closed.  By hypothesis (\ref{part-stable-reduction}),
$\tilde h(0) \in f^{-1}(\cV^+) \cap \cW^{-}$.  Thus, $f$ is stabilizer preserving at $\tilde{h}(0)$ and $f(\tilde{h}(0))$ is closed; therefore $\tilde{h}(0) \in \cU$. But this is a contradiction since we cannot have $\tilde h(\Spec K) \notin \cU$ and $\tilde h(0) \in \cU$.  Thus, $w \notin \bar{\cZ}$ as desired.

We now claim that $f$ maps closed points to closed points and is stabilizer preserving at all closed points of $z \in \cW$.  If $z \notin \cW^-$, then $z \in f^{-1}(\cV^-)$ so that $f$ is stabilizer preserving at $z$ and $f(z) \in \cV^-$ is closed.
If $z \in \cW^-$, then $f(z) \in \cX^-$ is closed.  If $f(z) \in \cX$ is not closed, then hypothesis (\ref{part-closed2}) implies $f(z) \in \cV^+$.  Thus, $f$ is stabilizer preserving at $z$ and $f(z) \in \cX$ is closed.  

Therefore, there exists an affine, \'etale and surjective morphism $f: \cW \to \cX$ such that $f: \cW \to \cX$ and $f|_{\cW^-}: \cW^- \to \cX^-$ are stabilizer preserving at closed points and send closed points to closed points.
By condition (5), the closure $\overline{ \{ x \} }$ admits a good moduli space; it follows from Theorem \ref{theorem-global-existence} that there exists a good moduli space $\cX \to X$.  It remains to show that $X^- \to X$ is proper.  Consider the commutative cube:
$${\def\objectstyle{\scriptstyle}
\def\labelstyle{\scriptstyle}
\xymatrix@=20pt{
                &\cW^- \ar@{^(->}[rr] \ar[dd] \ar[dl]          &               & \cW \ar[dd] \ar[dl] \\
\cX^- \ar@{^(->}[rr] \ar[dd]&                       & \cX \ar[dd]	  & \\
                &W^- \ar[rr] \ar[dl]                &               & W\ar[dl] \\
X^-\ar[rr]  &                           &X      &
}}$$
By Corollary \ref{corollary-cartesian}, the left and right squares are cartesian.  Since the diagram
$$\xymatrix{
\cW^- \ar[r] \ar[d]		& \cX' \ar[d] \\
X^- \times_X W \ar[r]		& X'
}$$
is cartesian, by uniqueness of good moduli spaces, we have that $X^- \times_X W \cong W^-$.  Since $W^- \to W$ is proper, by \'etale descent, $X^- \to X$ is also proper.
\end{proof}

The following lemma will be used to show that if $\bar{\cM}_{g,n}(A_k)$ admits a good moduli space, then so does $\bar{\cM}_{g,n}(A_k^+)$.
\begin{prop} \label{prop-inductive-existence2}
Let $\cX$ be an algebraic stack of finite type over $k$.  Let $\cX^{+} \subseteq \cX$ be an open substack.  Suppose that:
\begin{enumerate1}
\item There exists a good moduli space $\phi: \cX \to X$.
\item For every closed point $x \in \cX$, there exists a cartesian diagram
\begin{equation} \label{diagram-inductive}
\xymatrix{
\cW \ar[d]^f        & \cW^+ \ar@{_(->}[l] \ar[d]^{f|_{\cW^+}}   \ar @{} [dl] |{\square} \\
\cX                 & \cX^+ \ar@{_(->}[l]
}
\end{equation}
where $f$ is affine and \'etale, and there exist good moduli spaces $\cW \to W$ and $\cW^+ \to W^+$ such that $W^+ \to W$ is proper.
\end{enumerate1}
Then there exists a good moduli space $\cX^+ \to X^+$ such that $X^+ \to X$ is proper.
\end{prop}

\begin{proof}
By hypothesis (1) and Corollary \ref{corollary-openness}, there exists an \'etale, affine and surjective morphism $f: \cW \to \cX$ and a cartesian diagram
$$\xymatrix{
\cW \ar[d]^f \ar[r]  \ar @{} [dr] |{\square}    & W \ar[d] \\
\cX \ar[r]                          &X
}$$
where $\cW \to W$ is a good moduli space.  Clearly, $f|_{\cW^+}$ is stabilizer preserving.  To check that $f|_{\cW^+}: \cW^+ \to \cX^+$ sends closed points to closed points, let $ w \in \cW^+$ be a closed point.  We claim that if $\overline{ \{w\}} \subseteq \cW$ denotes the closure of $w$ in $\cW$, then $f(\overline{ \{w\}}) \subseteq \cX$ is closed.  First, observe that $f(\overline{ \{w\}}) \subseteq \overline{ \{ f(w) \} }$.  If $f(\overline{ \{w\}}) \subseteq \cX$ is not closed, then there exists a specialization $f(w) \rightsquigarrow x_0$ to a closed point $x_0 \in \cX$ with $x_0 \notin f(\overline{ \{w\}})$.  If $w \rightsquigarrow w_0$ is a specialization to a closed point $w_0 \in \cW$, then since $f$ sends closed points to closed points, $f(w) \rightsquigarrow f(w_0)$ is a specialization to a closed point $f(w_0) \in \cX$.  However, since $\cX$ admits a good moduli space, any $k$-point admits a unique specialization to a closed point.  Therefore, $x_0 = f(w_0)$ contradicting $x_0 \notin f(\overline{ \{w\}})$.  Therefore, $f(\overline{ \{w\}}) \subseteq \cX$ is closed.  It follows that $f(w) \in \cX^+$ is closed. Indeed, if there were a specialization $f(w) \rightsquigarrow x_1$ in $\cX^+$, then the fact that $f(\overline{ \{w\}}) \subseteq \cX$ is closed would imply that this specialization lifts to a specialization  $w \rightsquigarrow w_1$ in $\cW$ with $w_1 \in \cW^+$, a contradiction.  We conclude that $f|_{\cW^+}: \cW^+ \to \cX^+$ is stabilizer preserving and sends closed points to closed points. 

We now check that for any point $x \in \cX^+$, the closed substack $\cZ^+ = \overline{ \{x\} } \subseteq \cX^+$ admits a good moduli space.   Let $\cZ = \overline{ \{x\} } \subseteq \cX$ and consider the diagram obtained by base changing Diagram (\ref{diagram-inductive}) by $\cZ \hookarr \cX$
\begin{equation} \label{diagram-packet}
\xymatrix{
\cQ \ar[d]^{f|_{\cQ}}    & \cQ^+ \ar@{_(->}[l] \ar[d]^{f|_{\cQ^+}}   \ar @{} [dl] |{\square} \\
\cZ                 & \cZ^+ \ar@{_(->}[l]
}
\end{equation}
Since $f|_{\cQ}$ is stabilizer preserving and sends closed points to closed points, by Corollary \ref{corollary-openness} there is cartesian diagram 
$$\xymatrix{
\cQ \ar[d]^f \ar[r]  \ar @{} [dr] |{\square}    & Q \ar[d] \\
\cZ \ar[r]                          & \Spec k
}$$
Where $\cQ \to Q$ and $\cZ \to \Spec k$ are good moduli spaces and $Q \to \Spec k$ is \'etale.  Therefore, $Q = \coprod_i \Spec k$ so may assume in Diagram (\ref{diagram-packet}) that $\cQ \to \cZ$ is an isomorphism.  It follows that $\cQ^+ \to \cZ^+$ is an isomorphism and that $\cZ^+$ admits a good moduli space.
  
By Theorem \ref{theorem-global-existence}, there exists a good moduli space $\cX^+ \to X^+$ and the argument in the proof of Proposition \ref{prop-inductive-existence1} implies that $X^+ \to X$ is proper.
\end{proof}

\subsection{Existence}
We will apply Propositions \ref{prop-inductive-existence1} and \ref{prop-inductive-existence2} to show that $\bar{\cM}_{g,n}(A_k)$ and $\bar{\cM}_{g,n}(A_k^+)$  admit proper good moduli spaces.  In order to apply the lemma, we first show that $\bar{\cH}_{g,n}(A_k), \bar{\cS}_{g,n}(A_k)$ admits a good moduli space. This is accomplished inductively by Lemmas \ref{lemma-tacnode-existence}, \ref{lemma-H-existence} and \ref{lemma-S-existence}.  We first make the observation that $\bar{\cM}_{g,n}(A_k)$ is a global quotient stack of a smooth scheme by a connected algebraic group, so that $\bar{\cM}_{g,n}(A_k)$, $\bar{\cM}_{g,n}(A_k^+)$, $\bar{\cH}_{g,n}(A_k)$ and $\bar{\cS}_{g,n}(A_k)$  each admit \'etale GIT presentations by Proposition \ref{prop-quotient}.

\begin{lem} \label{lemma-tacnode-existence}
Let $k = 3$ and $r \ge 1$.  Suppose that $\bar{\cM}_{g',n'}(A_k)$ admits a proper good moduli space for $g' < 2r-1$ and $n' \in \{2,4\}$. Then $\bar{\cH}_{2r-1,2}(A_k)$ admits a proper good moduli space.
\end{lem}

\begin{proof}  
By \cite[Proposition 3.9]{asw}, $\bar{\cH}_{1,2}(A_k)$ admits a good moduli space which is a point.  For $i=1, \ldots, r-1$, consider the gluing maps
$$\begin{aligned}
\Psi_i^{(0)}: \bar{\cM}_{2(r-i)-1,4}(A_k) \times \bar{\cH}_{2i-1,2}(A_k) &\to \bar{\cH}_{2r-1,2}(A_k)\\
\Psi_i^{(j)}: \bar{\cM}_{j,2}(A_k)  \times \bar{\cM}_{2(r-i) - j,2}(A_k) \times \bar{\cH}_{2i-1,2}(A_k)  &\to \bar{\cH}_{2r-1,2}(A_k)  & \quad \text{for $j = 1, \ldots, g-rm-1$} \\
\end{aligned}$$
whose image collectively maps onto the locus of curves admitting a $A_3$-chain of length less than $r$. 
Each morphism $\Psi_i^{(j)}$ is clearly representable and quasi-finite.  Furthermore, by the argument of \cite[Proposition 4.16]{asw}, $\Psi_i^{(j)}$ satisfies the valuative criterion for properness; it follows that $\Psi_i^{(j)}$ is finite. By Proposition \ref{prop-finite-existence}, the images of $\Psi_i^{(j)}$ admit proper good moduli spaces.

Let $\cZ_{r} \subseteq \bar{\cH}_{2r-1,2}(A_k)$ be the closed stack consisting of curves with $A_3$-chains of length $r$.  There is a unique closed point $z \in |\cZ_{r}|$ which is the nodal union of $2r-1$ monomial $H_{1,2}$-bridges; see \cite[Section 5.2]{asw}.   Recall that  $\bar{\cS}_{1,2}(A_k)$ is the moduli stack of families of $2$-pointed curves of arithmetic genus 1 obtained by gluing two smooth rational curves along an $A_3$-singularities and marking each of the rational curves (see \cite[Definition 3.3]{asw}).  

By the argument of \cite[Proposition 2.11]{smyth_elliptic2}, there exists a $\times_{i=1}^{r-1} \bar{\cS}_{1,2}(A_k)$-bundle 
$$\cC \to \times_{i=1}^r \bar{\cH}_{1,2}(A_k)$$
of $A_k$-stable curves inducing a finite morphism
$$\cC \to \bar{\cH}_{2r-1,2}(A_k)$$
whose image is $\cZ_r$; the moduli stack $\cC$ parameterizes the possible ways to glue $r$ curves in $\bar{\cH}_{1,2}(A_k)$ at $r-1$ $A_k$-singularities at their marked points.  Since $\cC$ is defined by the stacky projectivization of a vector bundle over $\bar{\cH}_{1,2}(A_k)$, $\cC$ admits a good moduli space.  Moreover,  $\bar{\cH}_{1,2}(A_k)$ admits a good moduli space.  Therefore, by Proposition \ref{prop-finite-existence}, $\cZ_r$ admits a good moduli space.

Since $\bar{\cH}_{2r+1,2}(A_k)$ is covered by $\cZ_{r}$ and the images of $\Psi_i^{(j)}$, we conclude by Proposition \ref{prop-finite-existence} that $\bar{\cH}_{2r+1,2}(A_k)$ admits a proper good moduli space.
\end{proof}

\begin{lem} \label{lemma-H-existence}
Suppose that $\bar{\cM}_{g',n'}(A_k)$ admits a proper good moduli space for $g' < g$ and all $n'$.  Then $\bar{\cH}_{g,n}(A_k)$ admits a proper good moduli space.
\end{lem}

\begin{proof}
Suppose that $k = 2m$.  By \cite[Proposition 3.9]{asw}, $\bar{\cH}_{m,1}(A_k)$ admits a good moduli space which is a point.
The gluing maps
$$\begin{aligned}
 \bar{\cM}_{g-m,n+1}(A_k) \times  \bar{\cH}_{m,1}(A_k)&\to \bar{\cH}_{g,n}(A_k) \\
\bar{\cH}_{m,1}(A_k)  \times \bar{\cH}_{m,1}(A_k) &\to \bar{\cH}_{2m,0}(A_k) \\
\end{aligned}$$
satisfy the valuative criterion for properness by \cite[Lemma 4.15]{asw}.  It follows that they are finite.  Moreover, their images cover $\bar{\cH}_{g,n}(A_k)$.  By the inductive hypothesis and Prop \ref{prop-finite-existence}, we conclude that $\bar{\cH}_{g,n}(A_k)$ admits a proper good moduli space.
 
Suppose that $k = 3$.
 For $r = 1, \ldots, \lfloor \frac{g+1}{2} \rfloor$, consider the gluing maps
 $$\begin{aligned}
 \bar{\cM}_{g-2r,n+2}(A_k) \times \bar{\cH}_{2r-1,2}(A_k) &\to \bar{\cH}_{g,n}(A_k) \\
  \bar{\cM}_{j,n_1}(A_k) \times   \bar{\cM}_{g-2r-j+1,n_2}(A_k) \times \bar{\cH}_{2r-1,2}(A_k) &\to \bar{\cH}_{g,n}(A_k)  \quad \text{for $j =0, \ldots, g-2r+1$ and $n_1+n_2=n+2$}
\end{aligned}$$
which maps onto the locus of curves admitting an $A_k$-chain of length $r$ where the two end points are nodes which are non-separating and separating, respectively.  Consider also the gluing map for $n \ge 1$ and $r = 1, \ldots, \lfloor \frac{g+1}{2} \rfloor$
$$\bar{\cM}_{g-2r+1,n}(A_k) \times \bar{\cH}_{2r-1,2}(A_k) \to \bar{\cH}_{g,n}(A_k) \quad \text{for $n \ge 1$} $$
 which maps onto the locus of curves admitting an $A_3$-chain of length $r$ where one of the end points is a marked point.  Finally, for the degenerate case, we have the gluing map
$$\bar{\cH}_{rm,2}(A_k) \to \bar{\cH}_{rm,0}(A_k)$$
 where the end points of the $A_3$-chain are identified in a single node.  It follows from  \cite[Lemma 4.15]{asw} that these morphisms are finite.  Since the union of their images cover $\bar{\cH}_{g,n}(A_k)$ so we may apply Lemma \ref{lemma-tacnode-existence} and Proposition \ref{prop-finite-existence} to conclude that $\bar{\cH}_{g,n}(A_k)$ admits a proper good moduli space.
\end{proof}

\begin{lem} \label{lemma-S-existence}
Suppose that $\bar{\cM}_{g',n'}(A_k)$ admits a proper good moduli space for $g' < g$ and all $n'$.  Then $\bar{\cS}_{g,n}(A_k)$ admits a proper good moduli space.
\end{lem}

\begin{proof}
Suppose $k=2m$, and let $(C, \{p_i\}) \in \bar{\cS}_{g,n}(A_k)$ be a closed point.  By \cite[Proposition 5.14]{asw}, $C$ is the union of a genus  $g-rm$ curve $K$ with $r$ $H_{m,1}$ tails at nodes $q_1, \ldots, q_r \in K$ and $(K, \{p_i\} \cup \{q_j\}) \in \bar{\cM}_{g-rm,n}(A_k)$ is a closed point.  We will show that $(C, \{p_i\})$ is in the image of a finite morphism $\cC \to \bar{\cS}_{g,n}(A_k)$ where $\cC$ admits a proper good moduli space.  By Proposition \ref{prop-finite-existence}, it follows then that $\bar{\cS}_{g,n}(A_k)$ admits a proper moduli space.

For simplicity, we assume that $r=1$.  Recall that  $\bar{\cS}_{m,1}(A_k)$ is the moduli stack of families of $1$-pointed curves of arithmetic genus $m$ obtained by gluing an $A_k$-singularity at a marked rational curve \cite[Definition 3.3]{asw}).  By the argument of \cite[Proposition 2.11]{smyth_elliptic2}, there exists a stacky projective bundle
$$\cC \to \bar{\cM}_{g-rm,1}(A_k)$$
of $A_k$-stable curves inducing a finite morphism
$$\cC \to \bar{\cS}_{2r-1,2}(A_k)$$
whose image contains $(C, \{p_i\})$.   Indeed, in the case $k=2$, $\cC=\bar{\cM}_{g-rm,1}(A_k)$, and in the case $k=4$, $\cC$ is the stacky proj of the vector bundle $\psi_1^{3} \oplus \cO$ on $\bar{\cM}_{g-rm,1}(A_k)$; either way, $\cC$ admits a good moduli space.

Suppose $k=3$, and let $(C, \{p_i\}) \in \bar{\cS}_{g,n}(A_k)$ be a closed point.  For simplicity, let us assume that the $A_k$-chains of $C$ are non-disconnecting and attached at nodes (the other cases are essentially identical). By \cite[Proposition 5.14]{asw}, $C$ is the union of $k$ $A_3$-chains of length $r_1, \ldots, r_k$ and a genus $g-r-k$ core in $\bar{\cM}_{g-r-k,n+2k}(A_3^-)$, where $r=\sum_ir_i$. Note that $C$ is in the image of a finite morphism $\cC \to \bar{\cM}_{g,n}(A_k)$ where $\cC$ is a stacky projective bundle over $\times_{i=1}^{r} \bar{\cS}_{1,2}(A_k) \times \bar{\cM}_{g-r-k, n+2k}(A_k)$. Indeed $\cC$ is the fiber product of the $\PP^1$-bundles parametrizing one-dimensional subspaces of the direct sum of the cotangent spaces of the attaching points of the $r$ tacnodes. Since $\cC$ admits a good moduli space, we may apply Proposition \ref{prop-finite-existence} to conclude that $\bar{\cS}_{g,n}(A_k)$ admits a good moduli space.  We leave details to the reader.
\end{proof}

\begin{lem} \label{lemma-stable-reduction}
For every diagram
$$
\xymatrix{
 \Spec K \ar[r] \ar[d]  & \bar{\cM}_{g,n}(A_k^{-}) \ar@{^(->}[d] \\
 \Spec R \ar[r]^h \ar@{-->}[ur]^{\tilde h}  & \bar{\cM}_{g,n}(A_k)
}$$
with $R$ a valuation ring with fraction field $K$ such that $h(0)$ is a closed point where $0 \in \Spec R$ is the closed point and lift $\tilde h: \Spec R \to \bar{\cM}_{g,n}(A_k^{-})$ extending $\Spec K \to \bar{\cM}_{g,n}(A_k^-)$, then $\tilde h(0) \in \bar \cH_{g,n}(A_k)$.
\end{lem}

\begin{proof} This follows from analyzing the actions on the deformation space of $h(0)$ as in \cite[Section 6]{asw}.
\end{proof}

\begin{thm}  \label{theorem-main}
For $k=2,3$ and $4$, there exist good moduli spaces $\bar{\cM}_{g,n}(A_k) \to \bar{M}_{g,n}(A_k)$ and $\bar{\cM}_{g,n}(A_k^+) \to \bar{M}_{g,n}(A_k^+)$ inducing a commutative diagram
$$\xymatrix{
\bar{\cM}_{g,n}(A_k^{-}) \ar@{^(->}[r] \ar[d]   &\bar{\cM}_{g,n}(A_k) \ar[d]        & \bar{\cM}_{g,n}(A_k^{+}) \ar@{_(->}[l] \ar[d]\\
\bar{M}_{g,n}(A_k^{-})\ar[r]                & \bar{M}_{g,n}(A_k)                & \bar{M}_{g,n}(A_k^{+})\ar[l]
}$$
with $\bar{M}_{g,n}(A_k)$ and $\bar{M}_{g,n}(A_k^+)$ proper over $k$.
\end{thm}

\begin{proof}  Assume that the result holds for $k' < k$ so that there is a good moduli space $\bar{\cM}_{g,n}(A_k^-) \to M_{g,n}(A_k^-)$.  If $k$ is odd and $2g+1 < k$ or $k$ is even and $2g < k$, then $\bar{\cM}_{g,n}(A_k) = \bar{\cM}_{g,n}(A_{k-1}^+)$.  By induction, we may assume that $\bar{\cM}_{g',n'}(A_k)$ admits a proper good moduli space for $g'<g$ and all $n'$. 
We check the conditions of Proposition \ref{prop-inductive-existence1} where $\cX = \bar{\cM}_{g,n}(A_k)$, $\cX^- = \bar{\cM}_{g,n}(A_k^-)$,  $\cV^- =\bar{\cS}_{g,n}(A_k) $ and  $\cV^+ =\bar{\cH}_{g,n}(A_k)$.  For condition (1), $ \bar{\cM}_{g,n}(A_k^-)$ admits a good moduli space by the inductive hypothesis and $\bar{\cS}_{g,n}(A_k)$, $\bar{\cH}_{g,n}(A_k)$ admit good moduli spaces by Lemmas \ref{lemma-H-existence} and \ref{lemma-S-existence}.  Conditions (2) and (3) follow from the classisification of closed points in \cite[Section 5]{asw}.  Condition (4) follows from analyzing the actions on the deformation space of $h(0)$ as in \cite[Section 6]{asw}.  For condition (5), if $x \in \bar{\cM}_{g,n}(A_k)$ is any $k$-point, then the closure $\overline{ \{ x \} } \subseteq \bar{\cM}_{g,n}(A_k)$ is contained in $\bar{\cM}_{g,n}(A_k^-)$, $\bar{\cH}_{g,n}(A_k)$ or $\bar{\cS}_{g,n}(A_k)$ and therefore $\overline{ \{ x \} }$ admits a good moduli space since each of the substacks $\bar{\cM}_{g,n}(A_k^-)$, $\bar{\cH}_{g,n}(A_k)$ or $\bar{\cS}_{g,n}(A_k)$ do.  Condition (6) follows \cite[Theorem 8.3]{asw}.  Therefore there is a good moduli space $\bar{\cM}_{g,n}(A_k) \to M_{g,n}(A_k)$ with $\bar{M}_{g,n}(A_k^-) \to \bar{M}_{g,n}(A_k)$ proper.

We may now apply Proposition \ref{prop-inductive-existence2}, by using \cite[Theorem 8.3]{asw} again to verify condition (2) to conclude that there exists a good moduli space $\bar{\cM}_{g,n}(A_k^+) \to \bar{M}_{g,n}(A_k^+)$ with $\bar{M}_{g,n}(A_k^+) \to \bar{M}_{g,n}(A_k)$ proper.  
\end{proof}

\appendix
\section{}

In this appendix, we give examples of algebraic stacks including moduli stacks of curves which fail to have a good moduli space owing to a failure of conditions (1a), (1b), and (2) of Theorem B.  Note that there is an obviously necessary topological condition for a stack to admit a good moduli space, namely that every $k$-point has a unique isotrivial specialization to a closed point, and each of our examples satisfies this condition. The purpose of these examples is to illustrate the more subtle kinds of stacky behavior that can obstruct the existence of good moduli spaces.

\subsection*{Failure of condition (1a)}

\begin{example}  Let $\cX = [X/\ZZ_2]$ be the quotient stack where $X$ is the non-separated affine line and $\ZZ_2$ acts on $X$ by swapping the origins and fixing all other points.  The algebraic stack clearly satisfies condition (1b) and (2).  Then there is an \'etale, affine morphism $\AA^1 \to \cX$ which is stabilizer preserving at the origin but is not stabilizer preserving in an open neighborhood.  The algebraic stack $\cX$ does not admit a good moduli space.
\end{example}

While the above example may appear entirely pathological, we now provide two natural moduli stacks similar to this example.

\begin{example}
Consider the Deligne-Mumford locus $\cX \subseteq [\Sym^4 \PP^1 / \PGL_2]$ of unordered tuples $(p_1, p_2, p_3, p_4)$ where at least three points are distinct.  Consider the family  $(0, 1, \lambda, \infty)$ with $\lambda \in \PP^1$.  When $\lambda \notin \{ 0,1, \infty\}$, $\Aut(0,1, \lambda, \infty) \cong \ZZ/2\ZZ \times \ZZ/2\ZZ$; indeed, if $\sigma \in \PGL_2$ is the unique element such that $\sigma(0)=\infty$, $\sigma(\infty)=0$ and $\sigma(1) = \lambda$, then $\sigma([x,y]) = [y,\lambda x]$ so that $\sigma(\lambda)=1$ and therefore $\sigma \in \Aut(0,1, \lambda, \infty)$.  Similarly, there an element $\tau$ which acts via $0 \stackrel{\tau}{\leftrightarrow} 1$, $\lambda \stackrel{\tau}{\leftrightarrow} \infty$ and an element $\alpha$ which acts via $0 \stackrel{\alpha}{\leftrightarrow} \lambda$, $1 \stackrel{\tau}{\leftrightarrow} \infty$.  However, if $\lambda \in \{0,1,\infty\}$, $\Aut(0,1, \lambda, \infty) \cong \ZZ/2\ZZ$.

Therefore if $x=(0,1,\infty, \infty)$, any \'etale morphism $f: [\Spec A/ \ZZ_2] \to \cX$, where $\Spec A$ is a $\ZZ_2$-equivariant algebraization of the deformation space of $x$, will be stabilizer preserving at $x$ but not in any open neighborhood.  This failure of condition (1a) here is due to the fact that automorphisms of the generic fiber to not extend to the special fiber.  The algebraic stack $\cX$ does not admit a good moduli space but we note that if one enlarges the stack $\cX \subseteq [(\Sym^4 \PP^1)^{\ss} / \PGL_2]$ to include the point $(0,0,\infty,\infty)$, there does exist a good moduli space.

\end{example}

\begin{example}
Let $\cU_{2}$ be the stack of all reduced, connected curves of genus 2, and let $[C] \in \cU_{2}$ denote a cuspidal curve whose pointed normalization is a generic 1-pointed smooth elliptic curve $(E,p)$. We will show that any Deligne-Mumford open neighborhood $\cM \subset \cU_{2}$ of $[C]$ is non-separated and fails to satisfy condition (1a).

Note that $\Aut(C)=\Aut(E,p)=\ZZ/2\ZZ$. Thus, to show that no \'etale neighborhood 
$$[\Def(C)/\Aut(C)] \rightarrow \cM$$ can be stabilizer preserving, it is sufficient to exhibit a family $\cC \rightarrow \Delta$ whose special fiber is $C$, and whose generic fiber has automorphism group $\ZZ/2\ZZ \times \ZZ/2\ZZ$. To do this, let $C'$ be the curve obtained by nodally gluing two identical copies of $(E,p)$ along their respective marked points. Then $C'$ admits an involution swapping the two components, and a corresponding degree 2 map $C' \rightarrow E$ ramified over the single point $p$. We may smooth $C'$ to a family $\cC' \rightarrow \Delta$ of smooth double covers of $E$, simply by separating the ramification points. By \cite[Lemma 2.12]{smyth_elliptic1}, there exists a birational contraction 
$\cC' \rightarrow \cC$ contracting one of the two copies of $E$ in the central fiber to a cusp. The family $\cC \rightarrow \Delta$ now has the desired properties; the generic fiber has both a hyperelliptic and bielliptic involution while the central fiber is $C$.
\end{example}


\subsection*{Failure of condition (1b)}

\begin{example}
Let $\cX = [\AA^2 \setminus 0 / \GG_m]$ where $\GG_m$ acts via $t \cdot (x,y) = (x,ty)$.  Let $\cU = \{y \neq 0 \}  = [\Spec k[x,y]_y / \GG_m] \subseteq \cX$.   Observe that the point $(0,1)$ is closed in $\cU$ and $\cX$.  Then the open immersion $f: \cU \to \cX$ has the property that $f(0,1) \in \cX$ is closed but for $x \neq 0$, $(x,1) \in \cU$ is closed but $f(x,1) \in \cX$ is not closed.  There is no open neighborhood of $(0,1)$ in $\cX$ which admits a good moduli space.
\end{example}

\begin{example}
Let $\cM=\bar{\cM_{g}} \cup \cM^1 \cup \cM^2$, where $\cM^1$ consists of all curves of arithmetic genus $g$ with a single cusp and smooth normalization, and $\cM^2$ consist of all curves of the form $D \cup E_0$, where $D$ is a smooth curve of genus $g-1$ and $E_0$ is a rational cuspidal curve attached to $C$ nodally. 

We observe that $\cM$ has the following property: If $C=D \cup E$, where $D$ is a curve of genus g-1 and $E$ is a curve an elliptic tail, then $[C] \in \cM$ is a closed point if and only if $D$ is singular. Indeed, if $D$ is smooth, then $C$ admits an isotrivial specialization to $D \cup E_0$, where $E_0$ is a rational cuspidal tail.

Now consider any curve of the form $C=D \cup E$ where $D$ is a singular curve of genus $g-1$ and $E$ is a smooth elliptic tail, and, for simplicity, assume that $D$ has no automorphisms. We claim that there is no \'etale neighborhood of the form $[\Def(C)/\Aut(C)] \rightarrow \cM$, which sends closed points to closed points. Indeed, curves of the form $D' \cup E$ where $D'$ is smooth will appear in any such neighborhood and will obviously be closed in $[\Def(C)/\Aut(C)]$ (since this is a Deligne-Mumford stack), but are not closed in $\cM$.
\end{example}

\subsection*{Failure of condition (2)}

\begin{example} \label{example-nodal-cubic}
Let  $\cX = [C/\GG_m]$ where $C$ is the nodal cubic curve with the $\GG_m$-action given by multiplication.  Observe that $\cX$ is an algebraic stack with two points--one open and one closed.  But $\cX$ does not admit a good moduli space; if it did, $\cX$ would necessarily be cohomollogically affine and consequently $C$ would be affine, a contradiction.  However, there is an \'etale and affine morphism (but not finite) morphism $\cW=[\Spec (k[x,y]/xy)  / \GG_m] \to \cX$ where $\GG_m = \Spec k[t]_t$ acts on $\Spec k[x,y]/xy$ via $t \cdot (x,y) = (tx, t^{-1}y)$ which is stabilizer preserving and sends closed points to closed points;  however, the two projections $\cW \times_{\cX} \cW \rrarrows \cW$ do not send closed points to closed points.

To realize this \'etale local presentation concretely, we may express $C=X/\ZZ_2$ where $X$ is the union of two $\PP^1$'s with coordinates $[x_1,y_1]$ and $[x_2,y_2]$ glued via nodes at $0_1=0_2$ and $\infty_1 = \infty_2$ by the action of $\ZZ/2\ZZ$ where $-1$ acts via $[x_1,y_1] \leftrightarrow [y_2, x_2]$.  There is a  $\GG_m$-action on $X$ given by $t \cdot [x_1,y_1] = [tx_1,y_1]$ and $t \cdot [x_2, y_2] = [x_1, ty_1]$ which descends to an action on $C$.  We therefore have a 2-cartesian diagram
$$\xymatrix{
\ZZ_2 \times [X/\GG_m] \ar[r]^{p_1} \ar[d]^{p_2}	& [X/\GG_m] \ar[d] \\
[X/\GG_m] \ar[r]							& [C/\GG_m]
}$$
Let $p = \infty_2 = 0_1, q = \infty_1 = 0_2 \in X$.  If we let $\cW = [ (X \setminus \{p\} ) / \GG_m]$ and consider $f: \cW \to \cX$, we obtain  a 2-cartesian diagram
 $$\xymatrix{
 [(X \setminus \{p\})/\GG_m] \coprod [(X \setminus \{p,q\}) / \GG_m] \ar[r]^{\qquad \qquad p_1} \ar[d]^{p_2}	& [(X \setminus \{p\})/\GG_m] \ar[d]^f \\
[(X \setminus \{p\})/\GG_m] \ar[r]^f							& [C/\GG_m]
}$$

But  $[(X \setminus \{p,q\}) / \GG_m] \cong \Spec k \coprod \Spec k$ and the projections $p_1, p_2: \cW \times_{\cX} \cW \to \cW$ correspond to the inclusion of the two open points into $\cW$ which clearly doesn't send closed points to closed points.
\end{example}

\begin{example}[Condition (2)]
Let $C$ be the Deligne-Mumford semistable curve $D \cup E$, obtained by gluing a copy of $E:=\PP^1$ to a smooth genus $g-1$ curve $D$ at two points $p,q$. For simplicity, let us assume that $\Aut(D,p,q)=0$, so $\Aut(C)=\CC^*$. Let $\cM_{g}^{\ss,1}$ be the algebraic stack of Deligne-Mumford semistable curves $F$ where any rational subcurve connected to $F$ at only two points is smooth.

We will show that $\cM_{g}^{\ss,1}$ fails condition (2), i.e. the closed substack $\cZ:=\overline{ \{ [C] \} }$ fails to admit a good moduli space. It is easy to see that there is a unique isomorphism class of curves which isotrivially specializes to $C$, namely the nodal curve $C'$ obtained by gluing $D$ at $p$ and $q$. Thus, $\overline{ \{ [C] \} }$ has two points -- one open and one closed. We will show that $\overline{ \{ [C] \} }$ is isomorphic to the example given in Example \ref{example-nodal-cubic} of the quotient stack $[X/\GG_m]$ of the nodal cubic $X$ modulo $\GG_m$;  $[X/\GG_m]$ does not admit a good moduli space because $X$ is not affine.

To prove this, let us start by considering the constant family $D \times \PP^1$ with two constant sections corresponding to $p,q \in D$. Blowing up this family at $(p,0)$ and $(q, \infty)$, taking the strict transforms of the sections, and then identifying them nodally, we obtain a flat family of curves whose fibers are $C$ over $0$ and $\infty$ and $C'$ over every other point of $\PP^1$. The corresponding map $[\PP^1/\GG_m] \rightarrow \overline{ \{ [C] \} }$ is easily seen to factor through $[X/\GG_m]$, and the corresponding map $[X/\GG_m] \to \overline{ \{ [C] \} }$ is an isomorphism.
\end{example}


\bibliography{references}{}
\bibliographystyle{amsalpha}

\end{document}